\documentclass[journal,twoside,9pt]{IEEEtran}

\usepackage{graphicx,epsfig,subfigure,amsmath,amssymb,enumerate,subeqnarray}
\usepackage{srcltx,floatflt,wrapfig,psfrag,kotex}

\usepackage[sort]{cite}
\usepackage{dsfont,
           amssymb,
            enumerate,
            color,
            mathrsfs,
            subfigure,
             psfrag,
             epsfig,
            amsmath
            }

\usepackage{multicol}
\usepackage{epsfig,psfrag}

\newtheorem{notation}{Notation}

\newtheorem{theorem}{Theorem}

\newtheorem{remark}{Remark}

\newtheorem{problem}{Problem}
\newtheorem{lemma}{Lemma}
\newtheorem{definition}{Definition}

\allowdisplaybreaks

\begin{document}

%
\title{Relaxed Conditions for Parameterized Linear Matrix Inequality in the Form of Nested Fuzzy Summations}
%
%
%

\author{Do~Wan~Kim and Donghwan~Lee,~\IEEEmembership{Member,~IEEE}
 \thanks{D. W. Kim is with the Department
of Electrical Engineering, Hanbat National University, Deajeon,
Korea, e-mail: dowankim@hanbat.ac.kr.}

\thanks{D. Lee is with the Department
of Electrical Engineering, KAIST, Deajeon,
Korea, e-mail: donghwan@kaist.ac.kr.}

\thanks{\color{blue}This work has been submitted to IEEE Transactions on Systems, Man and Cybernetics: Systems for possible publications.}

}


\maketitle

\begin{abstract}
The aim of this study is to investigate less conservative conditions for parameterized linear matrix inequalities (PLMIs) that are formulated as nested fuzzy summations. Such PLMIs are commonly encountered in stability analysis and control design problems for Takagi-Sugeno (T-S) fuzzy systems. Utilizing the weighted inequality of arithmetic and geometric means (AM-GM inequality), we develop new, less conservative linear matrix inequalities for the PLMIs. This methodology enables us to efficiently handle the product of membership functions that have intersecting indices. Through empirical case studies, we demonstrate that our proposed conditions produce less conservative results compared to existing approaches in the literature.
\end{abstract}

\begin{IEEEkeywords}
T-S fuzzy system, relaxation, parameterized linear matrix inequality (PLMI), linear matrix inequality (LMI), stability
\end{IEEEkeywords}

\section{Introduction}

Parameterized linear matrix inequalities (PLMIs) play a crucial role in various disciplines, such as Takagi-Sugeno (T-S) fuzzy system analysis and control~\cite{Li2021,xie2022enhanced,wang2019fuzzy,ren2022static,jang2022intelligent,lu2021relaxed,Lee2014,lee2011approaches,Lee2014b}, and robust control designs. The objective is to identify matrix variables that satisfy the PLMIs over infinite-dimensional parameter spaces, commonly represented by matrix polytopes. In this paper, the challenge of searching for such matrix variables is referred to as the PLMI problem (PLMIP). Due to its infinite-dimensional nature, solving a PLMIP is generally computationally intractable, posing significant challenges.

To address this computational challenge, researchers have endeavored to establish finite-dimensional sufficient LMI conditions for PLMIPs using linear matrix inequalities (LMIs). For example,~\cite{Tuan2001,Kim2000,Teixeira2003,Fang2006,Tanaka2004,guerra2012} developed sufficient LMI conditions for PLMIPs expressed as double fuzzy sums (or two-fold nested fuzzy sums). Moreover, techniques based on Polya's theorem~\cite{young1934} have been utilized to extend these sufficient conditions to more general forms~\cite{Sala2007,Oliveira2007}.
Another avenue of research focuses on the sufficient LMI conditions set forth in~\cite{Tuan2001}, which are specifically designed for PLMIs expressed as two-fold nested fuzzy sums. The condition described in~\cite{Tuan2001} incorporates unique over-bounding techniques and employs a different set of vertices compared to other existing approaches.

Recently, a study by~\cite{Kim2023} introduced new, less conservative LMI conditions for PLMIPs formulated as two-fold nested fuzzy sums, extending the foundational concepts from~\cite{Tuan2001}. Specifically, the researchers proposed a novel sufficient LMI condition for PLMIPs based on Young's inequality. This newly developed condition encompasses the one presented in~\cite{Tuan2001} as a special case. However, the condition in~\cite{Kim2023} is limited to PLMIPs with two-fold nested fuzzy summations, constraining its applicability since many control design problems are formulated with more general nested fuzzy summations.

Motivated by this limitation, the main objective of this paper is to develop new, less conservative LMI conditions by generalizing the technique from~\cite{Kim2023} to accommodate general $q$-fold nested fuzzy summations. To achieve this, we employ the weighted inequality of arithmetic and geometric means (weighted AM-GM inequality) to efficiently handle the product of membership functions with intersecting indices. We empirically demonstrate through examples that the proposed conditions yield less conservative results compared to existing works based on Polya's theorem.
Finally, we note that the newly proposed LMI conditions distinguish themselves as a unique class, running parallel to those in~\cite{Sala2007,Oliveira2007}, which are based on Polya's theorem for addressing PLMIPs. Therefore, the proposed methods provide additional avenues for addressing PMLIPs, thereby broadening the range of solution techniques available.

\textit{Notations:} The notation $P\succ Q$ ($P\prec Q$) denotes that the matrix $P-Q$ is positive (negative) definite. In this context, $\mathbb{R}$ represents the set of real numbers, $\mathbb{R}^{n}$ denotes the $n$-dimensional Euclidean space, and $\mathbb{R}^{m\times n}$ refers to the field of real matrices with dimensions $m\times n$. Furthermore, $\mathbb{N}_{r}$ represents the positive integer set $\{1,2,\ldots,r\}$ where $r>1$. For two positive integers $a, b$, $\mathbb{N}_{[a,b]}$ denotes the set $\{a,a+1,\ldots,b\}$. To simplify notation, we will use $x$ instead of $x(t)$ for continuous-time signal vectors unless otherwise specified.

\section{Problem Formulation}

Consider the PLMI in the form of a $q$-fold nested fuzzy summations:
\begin{align}
\Phi(z)=\sum_{i_{1}=1}^{r}\sum_{i_{2}=1}^{r}\cdots\sum_{i_{q}=1}^{r}h_{i_{1}i_{2}\cdots i_{q}}(z)\Phi_{i_{1}i_{2}\cdots i_{q}}\prec0
\label{eq: PLMI}
\end{align}
where $h_{i_{1}i_{2}\cdots i_{q}}(z)=h_{i_{1}}(z)h_{i_{2}}(z)\times\cdots\times h_{i_{q}}(z)$. Here, $z\in\mathbb{R}^{p}$ represents the premise variable, $h_{i}(z):\mathscr{D}_{z}~\rightarrow \mathbb{R}_{[0,1]}$ satisfies $\sum_{i=1}^{r}h_{i}(z)=1$, $\mathscr{D}_{z}\subset \mathbb{R}^{p}$ denotes the domain of the premise variables, and $\Phi_{i_{1}i_{2}\cdots i_{q}}=\Phi_{i_{1}i_{2}\cdots i_{q}}^{T}\in\mathbb{R}^{n\times n},$ $(i_{1},i_{2},\cdots,i_{q})\in\mathbb{N}_{r}^{q}$ is a matrix that may have a linear dependency on the decision variables.

The PLMI~\eqref{eq: PLMI}, particularly when expressed as $2$-fold nested fuzzy summations, frequently arises in the context of controller design and stability analysis for linear uncertain models~\cite{Nguyen2018} and T-S fuzzy systems~\cite{Tuan2001,Tanaka2004}. To establish sufficient LMIs for negativeness, denoted as $\Phi_{i_{1}i_{2}}\prec0$ for all $(i_{1},i_{2})\in\mathbb{N}_{r}^{2}$, one can straightforwardly remove the terms $h_{i_{1}}h_{i_{2}}\in\mathbb{R}_{>0}$. Numerous efforts have been made to derive less conservative LMI conditions for the PLMI~\eqref{eq: PLMI} as $2$-fold nested fuzzy summations. Initially, relaxation techniques were introduced~\cite{Tuan2001,Tanaka2004}, leveraging properties of double fuzzy sums. For instance, these techniques exploit relationships like $\sum_{i_{1}=1}^{r}\sum_{i_{2}=1}^{r}h_{i_{1}}h_{i_{2}}\Phi_{i_{1}i_{2}}=\sum_{i_{1}=1}^{r}h_{i_{1}}^{2}\Phi_{i_{1}i_{1}}+\sum_{i=1}^{r}\sum_{i_{2}=1,i_{2}>i_{1}}^{r}h_{i_{1}}h_{i_{2}}(\Phi_{i_{1}i_{2}}+\Phi_{i_{2}i_{1}})$. These relaxation techniques~\cite{Tuan2001,Tanaka2004} lay the groundwork for subsequent research, including the introduction of slack variables~\cite{Kim2000,Teixeira2003}. As described in~\cite{guerra2006,guerra2012},  one of the well-known less conservative LMI conditions, which does not require the introduction of slack variables, is presented in the following lemma~\cite{Tuan2001}:
\begin{lemma}[Theorem~2.2 of \cite{Tuan2001}]\label{lemma: relaxation for dubled sum}
The PLMI~\eqref{eq: PLMI} in the form of $2$-fold nested fuzzy summations holds if
\begin{align*}
    \Phi_{i_{1}i_{1}} &\prec 0,\quad \forall i_{1}\in\mathbb{N}_{r}\\
    \frac{2}{r-1}\Phi_{i_{1}i_{1}}+\Phi_{i_{1}i_{2}}+\Phi_{i_{2}i_{1}}&\prec 0,~ \forall
    (i_{1},i_{2})\in\mathbb{N}_{r}^{2},~ i_{1}\neq i_{2}.
\end{align*}
\end{lemma}

Similarly, without the need for introducing slack variables, the authors have recently presented the following lemma, which definitely yields less conservative results than Lemma \ref{lemma: relaxation for dubled sum}:
\begin{lemma}[{Theorem~1 of~\cite{Kim2023}}]\label{lemma: 2-fold fuzzy summations}
The PLMI \eqref{eq: PLMI} holds in the form of $2$-fold nested fuzzy summations if
\begin{multline*}
    \Phi_{i_{1}i_{1}}
    +
    \frac{1}{2}\sum\limits_{i_{2} = 1,i_{2} < i_{1}}^r {\delta _{i_{2}} } (\Phi _{i_{1}i_{2}}  + \Phi _{i_{2}i_{1}} )
    \\
    + \frac{1}{2}\sum\limits_{i_{2} = 1,i_{2} > i_{1}}^r {\delta _{i_{2} - 1} } (\Phi _{i_{1}i_{2}}  + \Phi _{i_{2}i_{1}} ) \prec 0
\end{multline*}
for all $i_{1} \in {\mathbb N}_{r}$ and
$
(\delta _1 , \cdots ,\delta _{r - 1} ) \in \{0,1\}^{r-1}
$.
\end{lemma}

Unlike $2$-fold nested fuzzy summations, there are no relaxation results available for $q$-fold nested fuzzy summations without the use of slack variables, except for the following lemma (by replacing $\Phi_{i_{1}i_{2}}$ in Proposition 2 of \cite{Sala2007} by $\Phi_{i_{1}i_{2}\cdots i_{q}}$ ):
\begin{lemma}[\cite{Sala2007}]
The PLMI \eqref{eq: PLMI} in the form of $q$-fold nested fuzzy summations holds if
    \begin{align*}
        \sum_{(i_{1},i_{2},\cdots,i_{q})\in\mathcal{P}}
        \Phi_{i_{1}i_{2}\cdots i_{q}}
        \prec 0,\quad \forall(i_{1},i_{2},\cdots,i_{q})\in\mathbb{N}_{r}^{q}
    \end{align*}
where $\mathcal{P}$ is the set of all permutations of the indices $(i_{1},i_{2},\cdots,i_{q})$.
\end{lemma}

\begin{problem}
The primary objective of this paper is to formulate more generalized LMI conditions for PLMI \eqref{eq: PLMI}, extending the findings of Lemma~\ref{lemma: 2-fold fuzzy summations} from the context of $2$-fold fuzzy summations to encompass $q$-fold fuzzy summations.
\end{problem}

\section{Main Results}
To enhance clarity, we initially present the results for $3$-fold and $4$-fold fuzzy summations, followed by introducing a generalized result for $q$-fold case.

Before presenting our main results, we introduce the following lemma and definition, which will be used throughout the paper:
\begin{lemma}[{Weighted AM-GM inequality~\cite{gowers2010princeton}}]\label{lemma: am-gm}
 If  $c_{i}\in\mathbb{R}_{\geqslant 0}$, $i\in\mathbb{I}_{s}$ and  if
$\lambda_{i}\in\mathbb{R}_{> 0}$, $i\in\mathbb{I}_{s}$
such that $\sum_{i=1}^{s}\lambda_{i}=\lambda$,  then
\begin{align*}
    \sqrt[\lambda]{c_{1}^{\lambda_{1}}c_{2}^{\lambda_{2}}\times\cdots\times c_{s}^{\lambda_{s}}}
    \leqslant
    \frac{\lambda_{1}c_{1}+\lambda_{2}c_{2}+\cdots+\lambda_{s}c_{s}}{\lambda}
\end{align*}
in which equality holds if and only if all the $c_{i}$ with $\lambda_{i}$ are equal.
\end{lemma}

\begin{definition}\label{def: permutation}
Define the matrix-valued function
    $\mathscr{P}(\Phi_{i})
    =
    \sum_{i\in\mathbb{P}}\Phi_{i}$
for a given $\Phi_{i}$, where $\mathbb{P}$ denotes a set of permutations of the ordered $q$-tuple \(i=(i_{1},i_{2},\cdots,i_{q}) \in \mathbb{N}_{r}^{q}\). For example, when $i=(i_{1},i_{2})$, $\mathscr{P}(\Phi_{i_{1}i_{2}})=\sum_{i\in\{(i_{1},i_{2}),(i_{2},i_{1})\}}\Phi_{i}=\Phi_{i_{1}i_{2}}+\Phi_{i_{2}i_{1}}$.
\end{definition}
\begin{definition}\label{def: sigma}
Given the ordered $q$-tuple \(i\in \mathbb{N}_{r}^{q}\), we define a mapping \(\sigma(i)\) as follows:
\begin{multline*}
    \sigma(i):\left\{i\in\mathbb{N}_{r}^{q} : i_{a} \neq i_{b},\, a \in \mathbb{N}_{[2,q]},\, b \in \mathbb{N}_{a-1}\right\}
    \\
    \mapsto \left\{1,2,\cdots,\prod_{j=1}^{q-1}(r-j)\right\}.
\end{multline*}
To illustrate, let us consider the case of \(q=3\) and \(r=3\). When \(i_1=1\), we have two possibilities for \(i_2\) and \(i_3\) that are distinct from \(i_1\): i) \(i_2=2\) and \(i_3=3\), resulting in the mapping value \(\sigma((1,2,3))=1\);
ii) \(i_2=3\) and \(i_3=2\), resulting in the mapping value \(\sigma((1,3,2))=2\).
Thus, for \(i_1=1\), the domain of \(\sigma(i)\) is \(\{(1,2,3),(1,3,2)\}\), and the corresponding mapping values are \(\{1,2\}\). Similarly, for \(i_1=2\) and \(i_1=3\), we consider the possibilities of \(i_2\) and \(i_3\) that are distinct from \(i_1\), and determine the corresponding mapping values \(\sigma(i)\) for each case.
\end{definition}

\subsection{$3$-fold nested fuzzy summations}

The following lemma provides an equality that simplifies a $3$-fold nested fuzzy summation into a series of summations, particularly with respect to $\mathscr{P}(\Phi_{i_{1}i_{2}i_{3}})$ as defined in Definition \ref{def: permutation}:
\begin{lemma}\label{lemma: triple summation}
The following equality holds:
\begin{multline}
    \sum_{i_{1}=1}^{r}
    \sum_{i_{2}=1}^{r}
    \sum_{i_{3}=1}^{r}
     h_{i_{1}}h_{i_{2}}h_{i_{3}}\Phi_{i_{1}i_{2}i_{3}}
    \\
    =
    \sum_{i_{1}=1}^{r}
    h_{i_{1}}^{3}\Phi_{i_{1}i_{1}i_{1}}
    +
    \sum_{i_{1}=1}^{r}
    \sum_{\substack{i_{2}=1\\i_{2}\neq i_{1}}}^{r}
    h_{i_{1}}^{2}h_{i_{2}}
    \mathscr{P}(\Phi_{i_{1}i_{1}i_{2}})
    \\
    +
    \sum_{i_{1}=1}^{r}
    \sum_{\substack{i_{2}=1\\i_{2}\neq i_{1}}}^{r}
    \sum_{\substack{i_{3}=1\\i_{3}\neq i_{1}\\i_{3}\neq i_{2}}}^{r}
    h_{i_{1}}h_{i_{2}}h_{i_{3}}
    \frac{1}{3!}\mathscr{P}(\Phi_{i_{1}i_{2}i_{3}})
    \label{eq: triple summation}
\end{multline}
\end{lemma}
\begin{IEEEproof}
See Appendix \ref{app: lemma 3-fold}.
\end{IEEEproof}

The theorem presented below establishes sufficient LMI conditions for the PLMI formulated with $3$-fold nested fuzzy summations:
\begin{theorem}\label{th: LMI 3-dimensional}
The PLMI \eqref{eq: PLMI} in the form of $3$-fold nested fuzzy summations holds if
\begin{align}
\begin{aligned}
    &\Phi_{i_{1}i_{1}i_{1}}
    +
    \sum_{\begin{subarray}{c}i_{2}=1\\i_{2}\neq i_{1}\end{subarray}}^{r}
    \delta_{1,\sigma(i_{1},i_{2})}\frac{1}{3}
    (2\mathscr{P}(\Phi_{i_{1}i_{1}i_{2}})+\mathscr{P}(\Phi_{i_{2}i_{2}i_{1}}))
    \\&\quad
    + \sum_{\begin{subarray}{c}i_{2}=1\\i_{2}\neq i_{1}\end{subarray}}^{r}
    \sum_{\begin{subarray}{c}i_{3}=1\\i_{3}\neq i_{1}\\i_{3}\neq i_{2}\end{subarray}}^{r}
    \delta_{2,\sigma(i_{1},i_{2},i_{3})} \frac{1}{3!}
    \mathscr{P}(\Phi_{i_{1}i_{2}i_{3}})
    \prec 0
\end{aligned}
\label{eq: LMI for triple sums}
\end{align}
for all \(i_{1}\in\mathbb{I}_{r}\), \((\delta_{1,1}, \cdots,\delta_{1,r-1})\in\{0,1\}^{r-1}\), and \((\delta_{2,1}, \cdots,\delta_{2,(r-1)(r-2)})\in\{0,1\}^{(r-1)(r-2)}\), where a mapping $\sigma$ is given in Definition \ref{def: sigma}.
\end{theorem}

\begin{IEEEproof}
Utilizing Lemma \ref{lemma: triple summation}, we observe that
\begin{align*}
    \Phi(z)
    &\leqslant
    \sum_{i_{1}=1}^{r}
    h_{i_{1}}^{3}
    \Phi_{i_{1}i_{1}i_{1}}
     +\sum_{i_{1}=1}^{r}
    \sum_{\begin{subarray}{c}i_{2}=1\\i_{2}\neq i_{1}\end{subarray}}^{r}
    h_{i_{1}}^{2}h_{i_{2}}
    \max\left\{\mathscr{P}(\Phi_{i_{1}i_{1}i_{2}}),0\right\}
    \\
    &\quad +\sum_{i_{1}=1}^{r}
    \sum_{\begin{subarray}{c}i_{2}=1\\i_{2}\neq i_{1}\end{subarray}}^{r}
    \sum_{\begin{subarray}{c}i_{3}=1\\i_{3}\neq i_{1}\\i_{3}\neq i_{2}\end{subarray}}^{r}
    h_{i_{1}}h_{i_{2}}h_{i_{3}}
    \frac{1}{3!}
    \max\left\{\mathscr{P}(\Phi_{i_{1}i_{2}i_{3}}),0\right\}
\end{align*}
By utilizing the weighted AM-GM inequality from Lemma \ref{lemma: am-gm}, we see that
$h_{i_{1}}^{2}h_{i_{2}}\leqslant\frac{1}{3}(2h_{i_{1}}^{3}+h_{i_{2}}^{3})$ and $h_{i_{1}}h_{i_{2}}h_{i_{3}}\leqslant\frac{1}{3}(h_{i_{1}}^{3}+h_{i_{2}}^{3}+h_{i_{3}}^{3})$ in the foregoing inequality. From these properties, we deduce that
\begin{align*}
    \Phi(z)
    &\leqslant
    \sum_{i_{1}=1}^{r}
    h_{i_{1}}^{3}
    \Phi_{i_{1}i_{1}i_{1}}
    \\
    &\quad
    +\sum_{i_{1}=1}^{r}
    \sum_{\begin{subarray}{c}i_{2}=1\\i_{2}\neq i_{1}\end{subarray}}^{r}
    h_{i_{1}}^{3}
    \frac{1}{3}
    \max\left\{2\mathscr{P}(\Phi_{i_{1}i_{1}i_{2}})+\mathscr{P}(\Phi_{i_{2}i_{2}i_{1}}),0\right\} \\
    &\quad +\sum_{i_{1}=1}^{r}
    \sum_{\begin{subarray}{c}i_{2}=1\\i_{2}\neq i_{1}\end{subarray}}^{r}
    \sum_{\begin{subarray}{c}i_{3}=1\\i_{3}\neq i_{1}\\i_{3}\neq i_{2}\end{subarray}}^{r}
    h_{i_{1}}^{3}
    \frac{1}{3!}
    \max\left\{\mathscr{P}(\Phi_{i_{1}i_{2}i_{3}}),0\right\}
\end{align*}
where we have used
   $ h_{i_{2}}
    \mathscr{P}(\Phi_{i_{1}i_{1}i_{2}})
    = h_{i_{1}}
    \mathscr{P}(\Phi_{i_{2}i_{2}i_{1}})$ for all \((i_{1},i_{2})\in \{(i_{1},i_{2})\in\mathbb{N}_{r}^{2}:i_{2}\neq i_{1}\}\) and $h_{i_{1}}
    \mathscr{P}(\Phi_{i_{1}i_{2}i_{3}})
    = h_{i_{2}}
    \mathscr{P}(\Phi_{i_{2}i_{1}i_{3}})
    = h_{i_{3}}
    \mathscr{P}(\Phi_{i_{3}i_{2}i_{1}})$ for all \((i_{1},i_{2},i_{3})\in \{(i_{1},i_{2},i_{3})\in\mathbb{N}_{r}^{3}:i_{2}\neq i_{1},i_{3}\neq i_{1},i_{3}\neq i_{2}\}\). Thus, by employing \(\delta_{\sigma}^{1}\) and \(\delta_{\sigma}^{2}\) as discrete variables that take values of $0$ or $1$ instead of \(\max\), it can be deduced that the LMIs \eqref{eq: LMI for triple sums} imply \eqref{eq: PLMI}.
\end{IEEEproof}
\subsection{$4$-fold nested fuzzy summations}
The subsequent lemma presents an equation that simplifies a $4$-fold nested fuzzy summation into a sequence of summations, specifically focusing on the term $\mathscr{P}(\Phi_{i_{1}i_{2}i_{3}i_{4}})$ as defined in Definition \ref{def: permutation}:
\begin{lemma}\label{lemma: quadruple summation}
The following equality holds:
\begin{align}
\begin{aligned}
&\sum_{i_{1}=1}^{r}\sum_{i_{2}=1}^{r}\sum_{i_{3}=1}^{r}\sum_{i_{4}=1}^{r}
     h_{i_{1}}h_{i_{2}}h_{i_{3}}h_{i_{4}}
     \Phi_{i_{1}i_{2}i_{3}i_{4}} \\
&\quad =
\sum_{i_{1}=1}^{r}
h_{i_{1}}^{4}\Phi_{i_{1}i_{1}i_{1}i_{1}}
+
\sum_{i_{1}=1}^{r}
\sum_{\substack{i_{2}=1 \\ i_{2}\neq i_{1}}}^{r}
h_{i_{1}}^{3}h_{i_{2}}
\mathscr{P}(\Phi_{i_{1}i_{1}i_{1}i_{2}})
\\
&\quad
 +
\sum_{i_{1}=1}^{r}
\sum_{\substack{i_{2}=1 \\ i_{2}\neq i_{1}}}^{r}
h_{i_{1}}^{2}h_{i_{2}}^{2}
\frac{1}{2}\mathscr{P}(\Phi_{i_{1}i_{1}i_{2}i_{2}})\\
&\quad +
\sum_{i_{1}=1}^{r}
\sum_{\substack{i_{2}=1 \\ i_{2}\neq i_{1}}}^{r}
\sum_{\substack{i_{3}=1 \\ i_{3}\neq i_{1}\\i_{3}\neq i_{2}}}^{r}
h_{i_{1}}^{2}h_{i_{2}}h_{i_{3}}
\frac{1}{2}\mathscr{P}(\Phi_{i_{1}i_{1}i_{2}i_{3}})
\\
&\quad
+
\sum_{i_{1}=1}^{r}
\sum_{\substack{i_{2}=1 \\ i_{2}\neq i_{1}}}^{r}
\sum_{\substack{i_{3}=1 \\ i_{3}\neq i_{1}\\i_{3}\neq i_{2}}}^{r}
\sum_{\substack{i_{4}=1 \\ i_{4}\neq i_{1}\\i_{4}\neq i_{2}\\i_{4}\neq i_{3}}}^{r}
h_{i_{1}}h_{i_{2}}h_{i_{3}}h_{i_{4}}
\frac{1}{4!}\mathscr{P}(\Phi_{i_{1}i_{1}i_{2}i_{3}}).
\end{aligned}
\label{eq: quadruple sum}
\end{align}
\end{lemma}
\begin{IEEEproof}
See Appendix \ref{app: proof of lemma 4-fold}.
\end{IEEEproof}

The theorem below offers sufficient LMI conditions for the PLMI expressed as $4$-fold nested fuzzy summations:
\begin{theorem}\label{th: LMI 4-dimensional}
The PLMI \eqref{eq: PLMI} in the form of $4$-fold nested fuzzy summations holds if
\begin{align}
\begin{aligned}
    &\Phi_{i_{1}i_{1}i_{1}i_{1}}
    \\
    &\quad
    +
     \sum_{\begin{smallmatrix}i_{2}=1\\i_{2}\neq i_{1}\end{smallmatrix}}^{r}
     \delta_{1,1,\sigma(i_{1},i_{2})}
    \frac{1}{4}
    (
    3\mathscr{P}(\Phi_{i_{1}i_{1}i_{1}i_{2}})
    +\mathscr{P}(\Phi_{i_{2}i_{2}i_{2}i_{1}})
    )
       \\
   &\quad
    +
     \sum_{\begin{smallmatrix}i_{2}=1\\i_{2}\neq i_{1}\end{smallmatrix}}^{r}
     \delta_{1,2,\sigma(i_{1},i_{2})}
    \frac{1}{2}
    \mathscr{P}(\Phi_{i_{1}i_{1}i_{2}i_{2}})
   \\
   &\quad
   +
    \sum_{\begin{smallmatrix}
       i_{2}=1\\i_{2}\neq i_{1}
    \end{smallmatrix}}^{r}
    \sum_{\begin{smallmatrix}
        i_{3}=1\\i_{3}\neq i_{1}\\i_{3}\neq i_{2}
    \end{smallmatrix}}^{r}
    \delta_{2,1,\sigma(i_{1},i_{2},i_{3})}
    \frac{1}{8}
    (2\mathscr{P}(\Phi_{i_{1}i_{1}i_{2}i_{3}})
    \\
    &\quad
    +\mathscr{P}(\Phi_{i_{2}i_{2}i_{1}i_{3}})
    +\mathscr{P}(\Phi_{i_{3}i_{3}i_{2}i_{1}}))
    \\
    &\quad
    +
   \sum_{\begin{smallmatrix}
       i_{2}=1\\i_{2}\neq i_{1}
    \end{smallmatrix}}^{r}
    \sum_{\begin{smallmatrix}
        i_{3}=1\\i_{3}\neq i_{1}\\i_{3}\neq i_{2}
    \end{smallmatrix}}^{r}
    \sum_{\begin{smallmatrix}
        i_{4}=1\\i_{4}\neq i_{1}\\i_{4}\neq i_{2}\\i_{4}\neq i_{3}
    \end{smallmatrix}}^{r}
    \delta_{3,1,\sigma(i_{1},i_{2},i_{3},i_{4})}
    \frac{1}{4!}
    \mathscr{P}(\Phi_{i_{1}i_{2}i_{3}i_{4}})
    \prec0
\end{aligned}
\label{eq: LMI quadruple}
\end{align}
for  all $i_{1}\in\mathbb{N}_{r}$,
$(\delta_{1,1,1},\cdots,\delta_{1,1,r-1})
\in\{0,1\}^{r-1}$,
$(\delta_{1,2,1},\cdots,\delta_{1,2,r-1})
\in\{0,1\}^{r-1}$,
$(\delta_{2,1,1}, \cdots,$ $\delta_{2,1,(r-1)(r-2)})
\in\{0,1\}^{(r-1)(r-2)}$, and
$(\delta_{3,1,1}, \cdots,\delta_{3,1,(r-1)(r-2)(r-3)})
\in\{0,1\}^{(r-1)(r-2)(r-3)}$, where $\sigma$ is given in Definition \ref{def: sigma}.
\end{theorem}
\begin{IEEEproof}
Starting from Lemma \ref{lemma: quadruple summation}, we can simplify and rearrange $\Phi(z)$ in PLMI  to obtain the following upper bound:
\begin{equation*}
\resizebox{1\hsize}{!}{$
    \begin{split}
    \Phi(z)
    &\leqslant
     \sum_{i_{1}=1}^{r}
    h_{i_{1}}^{4}\Phi_{i_{1}i_{1}i_{1}i_{1}}
    +
    \sum_{i_{1}=1}^{r}
    \sum_{\begin{smallmatrix}i_{2}=1\\i_{2}\neq i_{1}\end{smallmatrix}}^{r}
    h_{i_{1}}^{3}h_{i_{2}}
    \max\left\{
    \mathscr{P}(\Phi_{i_{1}i_{1}i_{1}i_{2}})
    ,0\right\}
    \\
    &\quad
    +
    \sum_{i_{1}=1}^{r}
    \sum_{\begin{smallmatrix}i_{2}=1\\i_{2}\neq i_{1}\end{smallmatrix}}^{r}
    h_{i_{1}}^{2}h_{i_{2}}^{2}
    \frac{1}{2}
    \max\left\{
    \mathscr{P}(\Phi_{i_{1}i_{1}i_{2}i_{2}})
     ,0\right\}
    \\
    &\quad
    +
    \sum_{i_{1}=1}^{r}
      \sum_{\begin{smallmatrix}
       i_{2}=1\\i_{2}\neq i_{1}
    \end{smallmatrix}}^{r}
    \sum_{\begin{smallmatrix}
        i_{3}=1\\i_{3}\neq i_{1}\\i_{3}\neq i_{2}
    \end{smallmatrix}}^{r}
    h_{i_{1}}^{2}h_{i_{2}}h_{i_{3}}
    \frac{1}{2}
    \max\left\{
    \mathscr{P}(\Phi_{i_{1}i_{1}i_{2}i_{3}})
     ,0\right\}
    \\
    &\quad
    +
    \sum_{i_{1}=1}^{r}
    \sum_{\begin{smallmatrix}
       i_{2}=1\\i_{2}\neq i_{1}
    \end{smallmatrix}}^{r}
    \sum_{\begin{smallmatrix}
        i_{3}=1\\i_{3}\neq i_{1}\\i_{3}\neq i_{2}
    \end{smallmatrix}}^{r}
    \sum_{\begin{smallmatrix}
        i_{4}=1\\i_{4}\neq i_{1}\\i_{4}\neq i_{2}\\i_{4}\neq i_{3}
    \end{smallmatrix}}^{r}
    h_{i_{1}}h_{i_{2}}h_{i_{3}}h_{i_{4}}
    \frac{1}{4!}
    \max\left\{
    \mathscr{P}(\Phi_{i_{1}i_{2}i_{3}i_{4}})
     ,0\right\}
\end{split}$}
\end{equation*}
By applying the AM-GM inequality in Lemma \ref{lemma: am-gm}, we further proceed to
\begin{equation*}
\resizebox{1\hsize}{!}{$
    \begin{split}
\Phi(z)
    &\leqslant
     \sum_{i_{1}=1}^{r}
    h_{i_{1}}^{4}\Phi_{i_{1}i_{1}i_{1}i_{1}}
    +
    \sum_{i_{1}=1}^{r}
    \sum_{\begin{smallmatrix}i_{2}=1\\i_{2}\neq i_{1}\end{smallmatrix}}^{r}
    \frac{3h_{i_{1}}^{4}+h_{i_{2}}^{4}}{4}
    \max\left\{
    \mathscr{P}(\Phi_{i_{1}i_{1}i_{1}i_{2}})
    ,0\right\}
    \\
    &\quad
    +
    \sum_{i_{1}=1}^{r}
    \sum_{\begin{smallmatrix}i_{2}=1\\i_{2}\neq i_{1}\end{smallmatrix}}^{r}
    \frac{2h_{i_{1}}^{4}+2h_{i_{2}}^{4}}{4}
    \frac{1}{2}
    \max\left\{
    \mathscr{P}(\Phi_{i_{1}i_{1}i_{2}i_{2}})
     ,0\right\}
    \\
    &\quad
    +
    \sum_{i_{1}=1}^{r}
      \sum_{\begin{smallmatrix}
       i_{2}=1\\i_{2}\neq i_{1}
    \end{smallmatrix}}^{r}
    \sum_{\begin{smallmatrix}
        i_{3}=1\\i_{3}\neq i_{1}\\i_{3}\neq i_{2}
    \end{smallmatrix}}^{r}
    \frac{2h_{i_{1}}^{4}+h_{i_{2}}^{4}+h_{i_{3}}^{4}}{4}
    \frac{1}{2}
    \max\left\{
    \mathscr{P}(\Phi_{i_{1}i_{1}i_{2}i_{3}})
     ,0\right\}
    \\
    &\quad
    +
    \sum_{i_{1}=1}^{r}
    \sum_{\begin{smallmatrix}
       i_{2}=1\\i_{2}\neq i_{1}
    \end{smallmatrix}}^{r}
    \sum_{\begin{smallmatrix}
        i_{3}=1\\i_{3}\neq i_{1}\\i_{3}\neq i_{2}
    \end{smallmatrix}}^{r}
    \sum_{\begin{smallmatrix}
        i_{4}=1\\i_{4}\neq i_{1}\\i_{4}\neq i_{2}\\i_{4}\neq i_{3}
    \end{smallmatrix}}^{r}
    \frac{h_{i_{1}}^{4}+h_{i_{2}}^{4}+h_{i_{3}}^{4}+h_{i_{4}}^{4}}{4}
    \\
    &\quad\times
    \frac{1}{4!}
    \max\left\{
    \mathscr{P}(\Phi_{i_{1}i_{2}i_{3}i_{4}})
     ,0\right\}.
    \end{split}$}
\end{equation*}
In order to consolidate the above upper bound into the form of $\sum_{i_1}^{r}h_{i_1}^4$,  we make use of index changes and transform the terms involving $h_{i_{2}}$, $h_{i_{3}}$, and $h_{i_{4}}$ can be into terms involving $h_{i_{1}}$:
\begin{gather*}
    h_{i_{2}}^{4}\mathscr{P}(\Phi_{i_{1}i_{1}i_{1}i_{2}})|_{\substack{i_{2}\mapsto i_{1}\\i_{1}\mapsto i_{2}\\}}
    =
    h_{i_{1}}^{4}\mathscr{P}(\Phi_{i_{2}i_{2}i_{2}i_{1}})
    \\
    h_{i_{2}}^{4}\mathscr{P}(\Phi_{i_{1}i_{1}i_{2}i_{2}})|_{\substack{i_{2}\mapsto i_{1}\\i_{1}\mapsto i_{2}\\}}
    =
    h_{i_{1}}^{4}\mathscr{P}(\Phi_{i_{1}i_{1}i_{2}i_{2}})
    \\
    h_{i_{2}}^{4}
    \mathscr{P}(\Phi_{i_{1}i_{1}i_{2}i_{3}})|_{\substack{i_{2}\mapsto i_{1}\\i_{1}\mapsto i_{2}\\}}
    =
    h_{i_{1}}^{4}
    \mathscr{P}(\Phi_{i_{2}i_{2}i_{1}i_{3}})
    \\
    h_{i_{3}}^{4}
    \mathscr{P}(\Phi_{i_{1}i_{1}i_{2}i_{3}})|_{\substack{i_{3}\mapsto i_{1}\\i_{1}\mapsto i_{3}\\}}
    =
    h_{i_{1}}^{4}
    \mathscr{P}(\Phi_{i_{3}i_{3}i_{2}i_{1}})
     \\
     \begin{aligned}
          h_{i_{2}}^{4}
    \mathscr{P}(\Phi_{i_{1}i_{2}i_{3}i_{4}})|_{\substack{i_{2}\mapsto i_{1}\\i_{1}\mapsto i_{2}\\}}
    =
    h_{i_{3}}^{4}
    \mathscr{P}(\Phi_{i_{1}i_{2}i_{3}i_{4}})|_{\substack{i_{3}\mapsto i_{1}\\i_{1}\mapsto i_{3}\\}}
   \\
   =
    h_{i_{4}}^{4}
    \mathscr{P}(\Phi_{i_{1}i_{2}i_{3}i_{4}})|_{\substack{i_{4}\mapsto i_{1}\\i_{1}\mapsto i_{4}\\}}
    =
    h_{i_{1}}^{4}
    \mathscr{P}(\Phi_{i_{1}i_{2}i_{3}i_{4}})
     \end{aligned}
\end{gather*}
for all $(i_{1},i_{2},i_{3},i_{4})\in \{(i_{1},i_{2},i_{3},i_{4})\in\mathbb{N}_{r}^{4}:i_{2}\neq i_{1},i_{3}\neq i_{1},i_{3}\neq i_{2},i_{4}\neq i_{1},i_{4}\neq i_{2},i_{4}\neq i_{3}\}$. Therefore, by  applying these index changes, we arrive at
\begin{align*}
    \Phi(z)
    &\leqslant
    \sum_{i_{1}=1}^{r}
    h_{i_{1}}^{4}\Phi_{i_{1}i_{1}i_{1}i_{1}}
    \\
    &\quad
    +
    \sum_{i_{1}=1}^{r}
    \sum_{\begin{smallmatrix}i_{2}=1\\i_{2}\neq i_{1}\end{smallmatrix}}^{r}
    h_{i_{1}}^{4}
    \frac{1}{4}
    \max\left\{
    3\mathscr{P}(\Phi_{i_{1}i_{1}i_{1}i_{2}})
    +\mathscr{P}(\Phi_{i_{2}i_{2}i_{2}i_{1}})
    ,0\right\}
    \\
    &\quad
    +
    \sum_{i_{1}=1}^{r}
    \sum_{\begin{smallmatrix}i_{2}=1\\i_{2}\neq i_{1}\end{smallmatrix}}^{r}
    h_{i_{1}}^{4}
    \frac{1}{2}
    \max\left\{
    \mathscr{P}(\Phi_{i_{1}i_{1}i_{2}i_{2}})
     ,0\right\}
    \\
    &\quad
    +
    \sum_{i_{1}=1}^{r}
      \sum_{\begin{smallmatrix}
       i_{2}=1\\i_{2}\neq i_{1}
    \end{smallmatrix}}^{r}
    \sum_{\begin{smallmatrix}
        i_{3}=1\\i_{3}\neq i_{1}\\i_{3}\neq i_{2}
    \end{smallmatrix}}^{r}
    h_{i_{1}}^{4}\frac{1}{8}
    \max\left\{
    2\mathscr{P}(\Phi_{i_{1}i_{1}i_{2}i_{3}})
    \right.
    \\
    &\quad
    \left.
    +\mathscr{P}(\Phi_{i_{2}i_{2}i_{1}i_{3}})
    +\mathscr{P}(\Phi_{i_{3}i_{3}i_{2}i_{1}})
     ,0\right\}
    \\
    &\quad
    +
    \sum_{i_{1}=1}^{r}
    \sum_{\begin{smallmatrix}
       i_{2}=1\\i_{2}\neq i_{1}
    \end{smallmatrix}}^{r}
    \sum_{\begin{smallmatrix}
        i_{3}=1\\i_{3}\neq i_{1}\\i_{3}\neq i_{2}
    \end{smallmatrix}}^{r}
    \sum_{\begin{smallmatrix}
        i_{4}=1\\i_{4}\neq i_{1}\\i_{4}\neq i_{2}\\i_{4}\neq i_{3}
    \end{smallmatrix}}^{r}
    h_{i_{1}}^{4}
    \frac{1}{4!}
    \max\left\{
    \mathscr{P}(\Phi_{i_{1}i_{2}i_{3}i_{4}})
     ,0\right\}
\end{align*}
where the LMIs~\eqref{eq: LMI quadruple} is derived. Consequently, it can be concluded that \eqref{eq: PLMI} holds if LMIs~\eqref{eq: LMI quadruple} hold for $h_{i_{1}}^{4}\neq0$, where $\delta_{1,1,\sigma}, \delta_{1,2,\sigma}, \delta_{2,1,\sigma}, \delta_{3,1,\sigma}$ are discrete variables that takes $0$ or $1$, and they are used in place of the $\max$ function to simplify the expressions.
\end{IEEEproof}

\subsection{Generalization: $q$-fold nested fuzzy summation}
    The following notations, definitions and lemma are introduced to provide a concise representation for complex expressions involving multi-indices and nested summations, facilitating more streamlined and readable mathematical expressions:
\begin{notation}\label{notation: multi-index}
In the context of multi-indices, where $i\in\mathbb{N}_{ r}^{q}$ and $\lambda,\mu\in\mathbb{N}_{0}^{q}$, we introduce the following notations to simplify intricate mathematical expressions:
    $i^{\lambda}
    =\underbrace{i_{1}\cdots i_{1}}_{\lambda_{1}}\underbrace{i_{2}\cdots i_{2}}_{\lambda_{2}}\cdots\underbrace{i_{q}\cdots i_{q}}_{\lambda_{q}}$,
    $h_{i}^{\lambda}=h_{i_{1}}^{\lambda_{1}}h_{i_{2}}^{\lambda_{2}}\times\cdots \times h_{i_{q}}^{\lambda_{q}}$,
    $h_{i}=h_{i_{1}}h_{i_{2}}\times\cdots \times h_{i_{q}}$, and
    $\mu!=\mu_{1}!\mu_{2}!\times\cdots \times \mu_{q}!$.
For example, if $\lambda=(2,1,1,0)$, then $i^{\lambda}=i_{1}i_{1}i_{2}i_{3}$ and $h_{i}^{\lambda} =h_{i_{1}}^{2}h_{i_{2}}h_{i_{3}}$
\end{notation}
\begin{notation}\label{notation: multi-summation}
Let the following notation represent a nested summation over distinct indices:
\begin{gather*}
    \sum_{i_{1}=1}^{r}
    \sum_{
    \begin{smallmatrix}
    i_{j},i_{j}\neq i_{k}\\
    j\in\mathbb{N}_{[2,q]}\\
    k\in\mathbb{N}_{j-1}
    \end{smallmatrix}
    }
    =
     \sum_{i_{1}=1}^{r}
    \sum_{\begin{smallmatrix}i_{2}=1\\i_{2}\neq i_{1}\end{smallmatrix}}^{r}
    \cdots
    \sum_{
    \begin{smallmatrix}
    i_{q}=1\\
    i_{q}\neq i_{1}\\
    i_{q}\neq i_{2}\\
    \vdots\\
    i_{q}\neq i_{q-1}
    \end{smallmatrix}
    }^{r}.
\end{gather*}
\end{notation}

\begin{definition}\label{def: the set of ordered $q$-tuple}
 Let $\Lambda_{k}$, $k\in\mathbb{N}_{q}$ be the set of ordered $q$-tuple $\lambda=(\lambda_{1},\cdots,\lambda_{q})\in\mathbb{N}_{0}^{q}$ defined as $\Lambda_{k}
     =
     \{
     \lambda\in\mathbb{N}_{0}^{q}\mid
     \sum_{i=1}^{q}\lambda_{i}=q,~\sum_{i=1}^{q}[\lambda_{i}\neq0]=k,~ \lambda_1 \geq \lambda_2 \geq \ldots \geq \lambda_q
     \}$.
\end{definition}

\begin{definition}[{Stirling numbers~\cite{luo2011some}}]\label{def:stirling_numbers}
Consider  the set $\Lambda_{k}$, $k\in\mathbb{N}_{q}$ is defined in Definition \ref{def: the set of ordered $q$-tuple}. The Stirling numbers of the second kind, denoted as $s(q,k)$, enumerate the ways to partition a set of $q$ labeled objects into $k$ nonempty unlabeled subsets. These numbers can be computed using the formula
    $s(q,k) = \sum_{\lambda\in\Lambda_{k}}
    \binom{q}{\lambda_{1},\cdots,\lambda_{q}} \frac{1}{\mu(\lambda)!}$, where $ \mu_{j}(\lambda)=\sum_{k=1}^{q}[\lambda_{k}=j]$, $j\in\mathbb{N}_{q}$.
\end{definition}

\begin{lemma}\label{lemma:stirling_powers}
Consider the Stirling numbers of the second kind $s(q,k)$, as defined in Definition \ref{def:stirling_numbers}. Then, the powers of an indeterminate $r$ in terms of falling factorials can be expressed as $r^{q} = \sum_{k=0}^{q} s(q,k) \frac{r!}{(r-k)!}$.
\end{lemma}
\begin{definition}\label{def: general sigma function}
Let $\sigma=(\sigma_{1},\sigma_{2},\sigma_{3})$ be an ordered $3$-tuple. Here, $\sigma_{1}$, $\sigma_{2}$, and $\sigma_{3}$ are the mappings defined as
    \begin{itemize}
        \item The mapping $ \sigma_{1}(k):\mathbb{N}_{[2,q]}\mapsto\mathbb{N}_{q-1}$ is defined as $\sigma_{1}(k)=k-1$.
        \item The mapping $\sigma_{2}(k,\lambda):\Lambda_{k} \mapsto \left\{1,\cdots,|\Lambda_{k}|\right\}$ is defined  for each $\lambda\in \Lambda_{k}$, $k\in\mathbb{N}_{[2,q]}$,
where $\Lambda_{k}$ is defined in Definition \ref{def: the set of ordered $q$-tuple}.
        \item The mapping  $\sigma_{3}(k,i):\{i\in\mathbb{N}_{r}^{k} \mid i_{a} \neq i_{b},\, a \in \mathbb{N}_{[2,k]},\, b \in \mathbb{N}_{a-1}\} \mapsto \left\{1,2,\cdots,\prod_{j=1}^{k-1}(r-j)\right\}$ is defined for each $i\in\mathbb{N}_{r}^{k}$, $k\in\mathbb{N}_{[2,q]}$.
\end{itemize}
In summary, the ordered $3$-tuple $\sigma$
comprises one-to-one mappings from sets to sets, each serving the purpose of index labeling in different contexts.
\end{definition}

The following Lemma is the general version of Lemmas \ref{lemma: triple summation} and \ref{lemma: quadruple summation}:
\begin{lemma}\label{lemma: multidimensional summation}
The following equality holds:
\begin{align*}
    \sum_{i_{1}=1}^{r}
    \sum_{i_{2}=1}^{r}
    \cdots
    \sum_{i_{q}=1}^{r}
    h_{i}\Phi_{i}
    &=
    \sum_{\lambda\in\Lambda_{k},k\in\mathbb{N}_{q}}
    \frac{1}{\mu(\lambda)!}
      \sum_{i_{1}=1}^{r}
    \sum_{\substack{
    i_{a},
    i_{a}\neq i_{b}\\
    a\in\mathbb{N}_{[2,k]}\\
    b\in\mathbb{N}_{a-1}
    }
    }h_{i}^{\lambda}\mathscr{P}(\Phi_{i^{\lambda}}).
\end{align*}
\end{lemma}
\begin{IEEEproof}
See Appendix \ref{app: proof of Lemma q-fold}.
\end{IEEEproof}

The following theorem presents sufficient LMI conditions for the PLMI in the form of $q$-fold nested fuzzy summations:
\begin{theorem}\label{th: LMI q-dimensional}
The PLMI \eqref{eq: PLMI} holds if
\begin{multline}
    \Phi_{i_{1}\cdots i_{1}}
    \\
   +
    \sum_{\substack{\lambda\in\Lambda_{k}\\k\in\mathbb{N}_{[2,q]}}}
   \frac{1}{\mu(\lambda)!}
    \sum_{
    \begin{smallmatrix}
    i_{a},
    i_{a}\neq i_{b}\\
    a\in\mathbb{N}_{[2,k]}\\
    b\in\mathbb{N}_{a-1}
    \end{smallmatrix}
    }
    \delta_{\sigma}
      \frac{1}{q}
    \sum_{j=1}^{q}
    \lambda_{j}
    \left. \mathscr{P}(\Phi_{i^{\lambda}})\right|_{\begin{smallmatrix}
         i_{j}\mapsto i_{1}
         \\
         i_{1}\mapsto i_{j}
          \end{smallmatrix}}
    \prec0
    \label{eq: LMI q-dimensional}
\end{multline}
for  all  $i_{1}\in\mathbb{N}_{r}$ and $\delta_{\sigma}=\delta_{\sigma_{1},\sigma_{2},\sigma_{3}}\in\{0,1\}$, where $\sigma_{1}\in\mathbb{N}_{q-1}$, $\sigma_{2}\in \mathbb{N}_{|\Lambda_{k}|}$, and $\sigma_{3}\in\mathbb{N}_{\prod_{j=1}^{k-1}(r-j)}$ for each $k\in\mathbb{N}_{[2,q]}$. Here, $\Lambda_{k}$, $\mu(\lambda)$, and the subscript, 3-tuple $\sigma$ in $\delta_{\sigma}$ are given in Definitions \ref{def: the set of ordered $q$-tuple}, \ref{def:stirling_numbers}, and \ref{def: general sigma function}, respectively.
\end{theorem}
\begin{IEEEproof}
From Lemma \ref{lemma: multidimensional summation}, we see that
\begin{align*}
    \Phi
    &=
   \sum_{i_{1}=1}^{r}
    h_{i_{1}}^{q}
    \Phi_{i_{1}\cdots i_{1}}
   +
    \sum_{\substack{\lambda\in\Lambda_{k}\\k\in\mathbb{N}_{[2,q]}}}
   \frac{1}{\mu(\lambda)!}
      \sum_{i_{1}=1}^{r}
    \sum_{
    \begin{smallmatrix}
    i_{a},
    i_{a}\neq i_{b}\\
    a\in\mathbb{N}_{[2,k]}\\
    b\in\mathbb{N}_{a-1}
    \end{smallmatrix}
    }h_{i}^{\lambda}\mathscr{P}(\Phi_{i^{\lambda}})
      \\ &\leqslant
    \sum_{i_{1}=1}^{r}
    h_{i_{1}}^{q}
    \Phi_{i_{1}\cdots i_{1}}
    \\
    &\quad
   +
    \sum_{\substack{\lambda\in\Lambda_{k}\\k\in\mathbb{N}_{[2,q]}}}
   \frac{1}{\mu(\lambda)!}
      \sum_{i_{1}=1}^{r}
    \sum_{
    \begin{smallmatrix}
    i_{a},
    i_{a}\neq i_{b}\\
    a\in\mathbb{N}_{[2,k]}\\
    b\in\mathbb{N}_{a-1}
    \end{smallmatrix}
    }h_{i}^{\lambda}
   \max
    \left\{
    \mathscr{P}(\mathscr{P}(\Phi_{i^{\lambda}}),0
    \right\}
\end{align*}
where we used $\lambda=(q,0\cdots,0)\in\Lambda_{1}$.
By using AM-GM inequality in Lemma \ref{lemma: am-gm}, the term $h_{i}^{\lambda}
    \max
    \left\{
    \mathscr{P}(\mathscr{P}(\Phi_{i^{\lambda}}),0
    \right\}$ satisfies
\begin{equation*}
    \resizebox{1\hsize}{!}{$
    \begin{split}
    h_{i}^{\lambda}
    \max
    \left\{
    \mathscr{P}(\mathscr{P}(\Phi_{i^{\lambda}}),0
    \right\}
    &
    \leqslant
    \frac{1}{q}
    \sum_{j=1}^{q}
    \lambda_{j}h_{i_{j}}^{q}
     \max
    \left\{
    \mathscr{P}(\Phi_{i^{\lambda}}),0
    \right\}
     \\
     &=
    \frac{1}{q}
    \sum_{j=1}^{q}
    \lambda_{j}h_{i_{1}}^{q}
  \max
    \left\{
    \left. \mathscr{P}(\Phi_{i^{\lambda}})\right|_{\begin{smallmatrix}
         i_{j}\mapsto i_{1}
         \\
         i_{1}\mapsto i_{j}
          \end{smallmatrix}}
    ,0
    \right\}.
    \end{split}$}
\end{equation*}
From the foregoing inequality, we arrive at
\begin{align*}
    \Phi
    &\leqslant
      \sum_{i_{1}=1}^{r}
    h_{i_{1}}^{q}
    \Phi_{i_{1}\cdots i_{1}}
       +
    \sum_{\substack{\lambda\in\Lambda_{k}\\k\in\mathbb{N}_{[2,q]}}}
   \frac{1}{\mu(\lambda)!}
    \\
    &\quad
    \times
      \sum_{i_{1}=1}^{r}
    \sum_{
    \begin{smallmatrix}
    i_{a},
    i_{a}\neq i_{b}\\
    a\in\mathbb{N}_{[2,k]}\\
    b\in\mathbb{N}_{a-1}
    \end{smallmatrix}
    }
      \frac{1}{q}
    \sum_{j=1}^{q}
    \lambda_{j}h_{i_{1}}^{q}
  \max
    \left\{
    \left. \mathscr{P}(\Phi_{i^{\lambda}})\right|_{\begin{smallmatrix}
         i_{j}\mapsto i_{1}
         \\
         i_{1}\mapsto i_{j}
          \end{smallmatrix}}
    ,0
    \right\}.
\end{align*}
Therefore, employing the $\delta$ function, which assumes values of either 0 or 1,
to cover all potential cases of the $\mathrm{max}$ function, we obtain the resulting LMIs \eqref{eq: LMI q-dimensional} to ensure $\Phi\prec0$.
\end{IEEEproof}
\begin{remark}
Lemmas \ref{lemma: triple summation} and \ref{lemma: quadruple summation} correspond to specific instances of Lemma \ref{lemma: multidimensional summation} with $q=3$ and $q=4$, respectively. Likewise, Theorems \ref{th: LMI 3-dimensional} and \ref{th: LMI 4-dimensional} are particular cases of Theorem \ref{th: LMI q-dimensional} for $q=3$ and $q=4$, respectively. An illustrative example for $q=4$ can be found in Appendix \ref{app: example}.
\end{remark}

\section{An Example}
Consider the T--S fuzzy system $\dot{x}=\sum_{i=1}^{3}h_{i}(A_{i}x+B_{i}u)$ and $u=\sum_{i=1}^{3}h_{i}F_{i}Q^{-1}x$
consisting of three subsystems, each defined by a pair of matrices
\begin{align*}
A_{1}=&\begin{bmatrix}
    1.59&-7.29\\0.01&0
    \end{bmatrix},
    \quad
B_{1}=\begin{bmatrix}
    1\\0\end{bmatrix}\\
A_{2}=&\begin{bmatrix}
    0.02&-4.64\\0.35&0.21
    \end{bmatrix},
    \quad B_{2}=\begin{bmatrix}
    8\\0
    \end{bmatrix}\\
A_{3}=& \begin{bmatrix}
    -a&-4.33\\0&0
    \end{bmatrix}, \quad B_{3} =\begin{bmatrix}
    -b+6\\-1
    \end{bmatrix}
\end{align*}
taken from \cite{Fang2006}, where the parameters $a$ and $b$ are introduced for the purpose of assessing the feasibility of LMIs. The widely-accepted asymptotic stabilization condition for this fuzzy system, as presented in\cite{Tanaka2004}, is given by PLMI \eqref{eq: PLMI} with
    $\Phi_{i_{1}i_{2}}
    =
    (A_{i_{1}}Q+B_{i_{1}}F_{i_{2}})^{T}+A_{i_{1}}Q+B_{i_{1}}F_{i_{2}}
    $, $(i_{1},i_{2})\in \mathbb{N}_{3}^{2}$
in the form of fuzzy double sum. However, discussions within the multi-index approach \cite{Sala2007}, based on Polya's theorem, indicate the potential for obtaining less conservative results for the asymptotic stabilization condition. This can be accomplished by reformulating the stabilization condition in the form of nested fuzzy summations rather than the initially specified double fuzzy sum.

We investigate the feasibility of the LMIs outlined in Theorem \ref{th: LMI q-dimensional} under two distinct scenarios: one corresponding to $3$-fold nested fuzzy summations, $q=3$ (referred to as Theorem \ref{th: LMI 3-dimensional}), and the other to $4$-fold nested fuzzy summations, $q=4$ (referred to as Theorem \ref{th: LMI 4-dimensional}). This analysis spans the parameter space $(a,b)\in[0,10]\times[0,10]$.
In both cases, we enforce the condition $\Phi_{i_{1}i_{2}i_{3}}=\Phi_{i_{1}i_{2}}$ for $3$-fold nested fuzzy summations and $\Phi_{i_{1}i_{2}i_{3}i_{4}}=\Phi_{i_{1}i_{2}}$ for $4$-fold nested fuzzy summations. Additionally, for the sake of comparison, we assess the feasibility of LMIs (Polya conditions) as described in \cite[Proposition 2]{Sala2007} for both scenarios.

\begin{figure}
\centering
\includegraphics[width=0.4\textwidth]{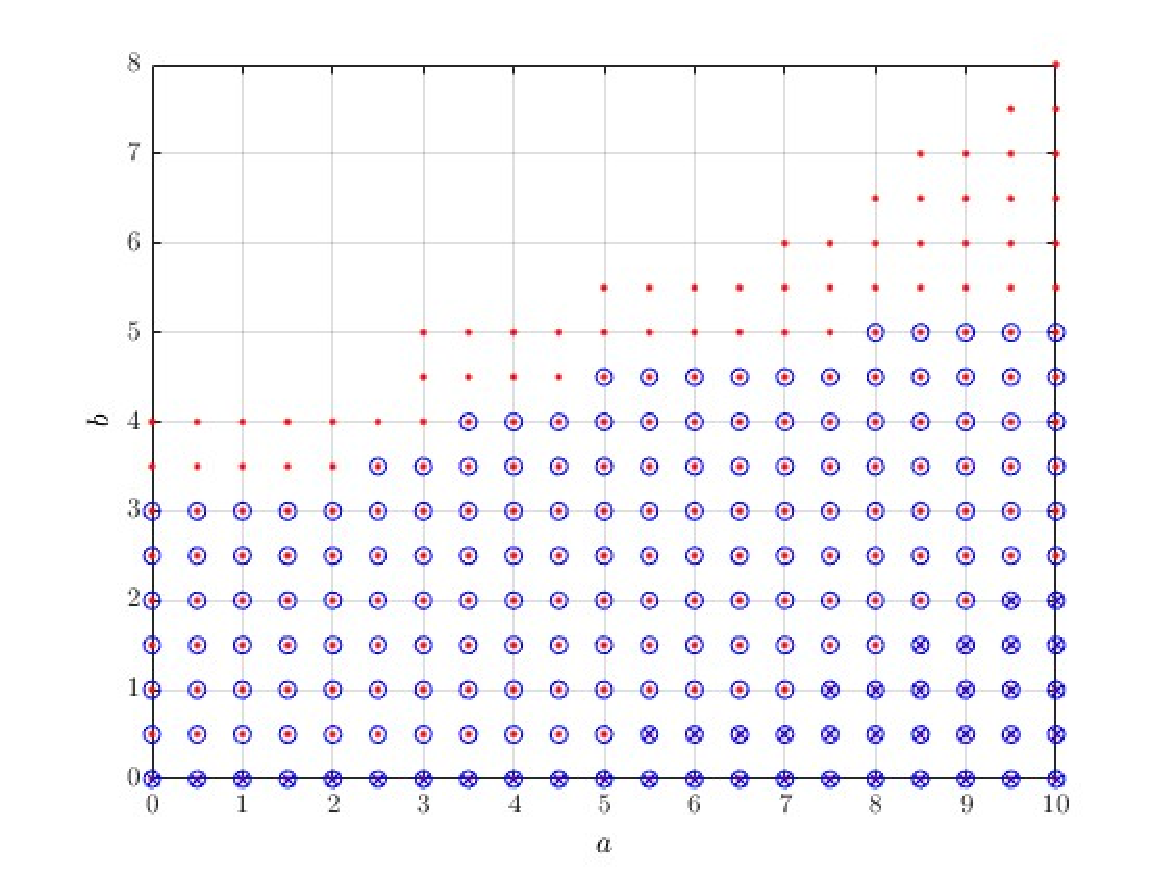}
\caption{Example~\label{ex:1}: Stabilizable region by Theorem~\ref{th: LMI q-dimensional} with $q=3$ (`$\circ$'), Theorem~\ref{th: LMI q-dimensional} with $q=4$ (`$\cdot$'), and \cite[Proposition 2]{Sala2007} (`$\times$') }\label{fig:1}
\end{figure}

Fig. \ref{fig:1} visually presents the regions of feasible points in the $(a,b)$ parameter space for both Theorem \ref{th: LMI q-dimensional} and \cite[Proposition 2]{Sala2007}. In this representation,  the marker `$\circ$' indicates the feasibility of the LMIs from Theorem \ref{th: LMI q-dimensional} concerning $3$-fold nested fuzzy summations case, while the marker `$\cdot$' signifies the feasibility of the LMIs for $4$-fold nested fuzzy summations case. Furthermore, the marker `$\times$' denotes the feasibility of the LMIs derived from \cite[Proposition 2]{Sala2007} for $4$-fold nested fuzzy summations case.

Remarkably, our feasibility assessments unveil that the LMIs proposed in \cite[Proposition 2]{Sala2007} for $3$-fold nested fuzzy summations scenario prove infeasible across the entire parameter space $(a,b)\in[0,10]\times[0,10]$. As illustrated in Fig. \ref{fig:1}, it becomes evident that Theorem~\ref{th: LMI q-dimensional} provides less conservative results in both $3$-fold and $4$-fold nested fuzzy summations scenarios compared to \cite[Proposition 2]{Sala2007} without the need for slack variables.

\section{Conclusions}

This paper introduces novel sufficient LMI conditions for PLMIs, formulated through general $q$-fold nested fuzzy summations, using the weighted AM-GM inequality. These conditions are derived without resorting to slack variables. Using experimental results, we have demonstrated that these conditions are less conservative compared to existing methodologies.

\appendices
\section{Proof of Lemma \ref{lemma: triple summation}}\label{app: lemma 3-fold}
We begin by expressing the summation involving triples as
$\sum_{i_{1}=1}^{r}
\sum_{i_{2}=1}^{r}
\sum_{i_{3}=1}^{r}
h_{i_{1}}h_{i_{2}}h_{i_{3}}\Phi_{i_{1}i_{2}i_{3}}
=
\sum_{i\in I^{3}}
h_{i_{1}}h_{i_{2}}h_{i_{3}}$ $\times \Phi_{i_{1}i_{2}i_{3}}$, where $ I^{3}$ represents the set of ordered triples $i=(i_{1}, i_{2}, i_{3})$ defined as $ I^{3}=\{i\in
    \mathbb{N}_{r}^{3}
    \}$. We categorize the possible combinations of indices $i_{1}$, $i_{2}$, and $i_{3}$ within $ I^{3}$ into three cases:
i) all three indices are equal;
ii) two of the indices are equal while the third one is distinct;
iii) all three indices are distinct.
In case i), there exists $\left(\begin{smallmatrix}
    3\\3
\end{smallmatrix}\right)(=1)$ combination ($i_{1}=i_{2}=i_{3}$). In case ii), there are $\left(\begin{smallmatrix}
    3\\2
\end{smallmatrix}\right)\left(\begin{smallmatrix}
    1\\1
\end{smallmatrix}\right)(=3)$ combinations ($i_{1}=i_{2}\neq i_{3}$, $i_{1}=i_{3}\neq i_{2}$, $i_{3}=i_{2}\neq i_{1}$). In case iii), there is $\left(\begin{smallmatrix}
    3\\1
\end{smallmatrix}\right)\left(\begin{smallmatrix}
    2\\1
\end{smallmatrix}\right)\left(\begin{smallmatrix}
    1\\1
\end{smallmatrix}\right)\frac{1}{3!}(=1)$ combination ($i_{1}\neq i_{2}\neq i_{3}$). For each of these cases, we define separate sets as follows: $ I_{1}^{3}=\{(i_{1},i_{2},i_{3})\in
    \mathbb{N}_{r}^{3}:i_{1}=i_{2}=i_{3}\}$, $I_{2}^{3}
     =
    \{(i_{1},i_{2},i_{3})\in
    \mathbb{N}_{r}^{3}:i_{1}=i_{2}\neq i_{3} \vee i_{1}=i_{3}\neq i_{2}
    \vee i_{2}=i_{3}\neq i_{1}\}$, and $I_{3}^{3}=
    \{(i_{1},i_{2},i_{3})\in
    \mathbb{N}_{r}^{3}:i_{1}\neq i_{2}\neq i_{3}\}$.

We observe that the original set $ I^{3}$ is the union of these three sets: $ I^{3}
    =
    \bigcup_{j=1}^{3}
     I_{j}^{3}$
since $\sum_{j=1}^{3}\left| I_{j}^{3}\right|=r+3r(r-1)+r(r-1)(r-2)=r^{3}=\left| I^{3}\right|$ and $\bigcup_{j=1}^{3}
     I_{j}^{3}\subseteq I^{3}$. This conclusion can be shown by contradiction. Suppose that $     I^{3}
    \neq
    \bigcup_{j=1}^{3}
     I_{j}^{3}$. Then, there exists $i\in  I^{3}$ such that $i\notin \bigcup_{j=1}^{3}
     I_{j}^{3}$. At this point, $\left| I^{3}\right|>\sum_{j=1}^{3}\left| I_{j}^{3}\right|$, but $\left| I^{3}\right|=\sum_{j=1}^{3}\left| I_{j}^{3}\right|$, which is a contradiction.

From this set equality and $\mathscr{P}$ in Definition \ref{def: permutation}, we can rewrite $3$-fold nested fuzzy summations as
\begin{equation*}
    \resizebox{1\hsize}{!}{
    $\begin{split}
    \sum_{i\in I^{3}}
    h_{i_{1}}h_{i_{2}}h_{i_{3}}\Phi_{i_{1}i_{2}i_{3}}
    &=
     \sum_{j=1}^{3}
     \sum_{i\in I_{j}^{3}}
     h_{i_{1}}h_{i_{2}}h_{i_{3}}\frac{1}{3!}\mathscr{P}(\Phi_{i_{1}i_{2}i_{3}})
    \\
    &
    =
    \sum_{i\in I_{1}^{3}}
    h_{i_{1}}^{3}\Phi_{i_{1}i_{1}i_{1}}
    +
    \sum_{i\in I_{2}^{3}}
    \left(
    h_{i_{1}}^{2}h_{i_{3}}\frac{2!}{3!}\mathscr{P}(\Phi_{i_{1}i_{1}i_{3}})
    \right.
    \\
    &\quad
    \left.
     +
    h_{i_{1}}^{2}h_{i_{2}}\frac{2!}{3!}\mathscr{P}(\Phi_{i_{1}i_{2}i_{1}})
    +
    h_{i_{1}}h_{i_{2}}^{2}\frac{2!}{3!}\mathscr{P}(\Phi_{i_{1}i_{2}i_{2}})
    \right)
    \\
    &\quad
    +
     \sum_{i\in I_{3}^{3}}
    h_{i_{1}}h_{i_{2}}h_{i_{3}}\frac{1}{3!}\mathscr{P}(\Phi_{i_{1}i_{2}i_{3}}).
    \end{split}$}
\end{equation*}
Hence, utilizing the properties $\sum_{i\in I_{2}^{3}}
    h_{i_{1}}^{2}h_{i_{3}}\mathscr{P}(\Phi_{i_{1}i_{1}i_{3}})
    =
     \sum_{i\in I_{2}^{3}}
    h_{i_{1}}^{2}h_{i_{2}}\mathscr{P}(\Phi_{i_{1}i_{1}i_{2}})$,
    $\sum_{i\in I_{2}^{3}}
    h_{i_{1}}h_{i_{2}}^{2}\mathscr{P}(\Phi_{i_{1}i_{2}i_{2}})
    =
     \sum_{i\in I_{2}^{3}}
    h_{i_{1}}$ $\times h_{i_{2}}^{2}\mathscr{P}(\Phi_{i_{2}i_{1}i_{1}})$, and
     $\sum_{i\in I_{2}^{3}}
     (
    \mathscr{P}(\Phi_{i_{1}i_{1}i_{2}})
    +
    \mathscr{P}(\Phi_{i_{1}i_{2}i_{1}})
    +
    \mathscr{P}(\Phi_{i_{2}i_{1}i_{1}})
    )
    =
    \sum_{i\in I_{2}^{3}}
    3
    \mathscr{P}(\Phi_{i_{1}i_{1}i_{2}})$,
we obtain the equation \eqref{eq: triple summation}.

\section{Proof of Lemma \ref{lemma: quadruple summation}}\label{app: proof of lemma 4-fold}
Similar to the proof process of Lemma \ref{lemma: triple summation}, we commence by expressing $4$-fold nested fuzzy summations as
$\sum_{i_{1}=1}^{r}
\sum_{i_{2}=1}^{r}
\sum_{i_{3}=1}^{r}
\sum_{i_{4}=1}^{r}
h_{i_{1}}h_{i_{2}}h_{i_{3}}h_{i_{4}}\Phi_{i_{1}i_{2}i_{3}i_{4}}
=
\sum_{i\in I^{4}}
h_{i_{1}}$ $\times h_{i_{2}}h_{i_{3}}h_{i_{4}}\Phi_{i_{1}i_{2}i_{3}i_{4}}$,
where $ I^{4}$ represents the set of ordered quadruples $i\in\mathbb{N}_{r}^{4}$. We classify the potential combinations of indices $i_{1}$, $i_{2}$, $i_{3}$, and $i_{4}$ within $ I^{4}$ into four distinct categories:
i) All four indices are equal;
ii) Three indices are equal while the fourth one is distinct, or there are two pairs that are identical within each pair but distinct between pairs;
iii) Two pairs of indices are equal while the remaining indices are distinct from each other;
iv) All four indices are distinct. In case i), there exists $\binom{4}{4}(=1)$ combination ($i_{1}=i_{2}=i_{3}=i_{4}$). In case ii), there are $\binom{4}{3}\binom{1}{1}(=4)$ combinations ($i_{1}=i_{2}=i_{3}\neq i_{4}$, $i_{1}=i_{2}=i_{4}\neq i_{3}$, $i_{1}=i_{3}=i_{4}\neq i_{2}$, $i_{2}=i_{3}=i_{4}\neq i_{1}$) or $\binom{4}{2}\binom{2}{2}\frac{1}{2!}(=3)$ combinations ($i_{1}=i_{2}\neq i_{3}= i_{4}$, $i_{1}=i_{3}\neq i_{2}=i_{4}$, $i_{1}=i_{4}\neq i_{2}=i_{3}$).
In case iii), there are $\binom{4}{2}\binom{2}{1}\binom{1}{1}\frac{1}{2!}(=6)$ combinations ($i_{1}=i_{2}\neq i_{3}\neq i_{4}$, $i_{1}=i_{3}\neq i_{2}\neq i_{4}$, $i_{1}=i_{4}\neq i_{2}\neq i_{3}$, $i_{2}=i_{3}\neq i_{1}\neq i_{4}$, $i_{2}=i_{4}\neq i_{1}\neq i_{3}$, $i_{3}=i_{4}\neq i_{1}\neq i_{2}$). In case iv), there is $\binom{4}{1}\binom{3}{1}\binom{2}{1}\binom{1}{1}\frac{1}{4!}(=1)$ combination ($i_{1}\neq i_{2}\neq i_{3}\neq i_{4}$).

To delineate each of these cases, we introduce distinct sets as follows: $I_{1}^{4}=\{i\in\mathbb{N}_{r}^{4}: i_{1}=i_{2}=i_{3}=i_{4}\}$, $I_{2}^{4}=I_{2,1}^{4}\cup I_{2,2}^{4}$,
$I_{3}^{4}:=\{i\in\mathbb{N}_{r}^{4}:i_{1}=i_{2}\neq i_{3}\neq i_{4}\vee i_{1}=i_{3}\neq i_{2}\neq i_{4}
    \vee  i_{1}=i_{4}\neq i_{2}\neq i_{3}
    \vee i_{2}=i_{3}\neq i_{1}\neq i_{4}
    \vee i_{2}=i_{4}\neq i_{1}\neq i_{3} \vee i_{3}=i_{4}\neq i_{1}\neq i_{2}\}$, and $I_{4}^{4}=\{i\in\mathbb{N}_{r}^{4}:i_{1}\neq i_{2}\neq i_{3}\neq i_{4}\}$,
where $I_{2,1}^{4}
    =
     \{i\in\mathbb{N}_{r}^{4}:
    i_{1}=i_{2}=i_{3}\neq i_{4} \vee i_{1}=i_{2}=i_{4}\neq i_{3}
    \vee i_{1}=i_{3}=i_{4}\neq i_{2}\vee i_{2}=i_{3}=i_{4}\neq i_{1}
    \}$ and $I_{2,2}^{4}
    =
    \{i\in\mathbb{N}_{r}^{4}:
    i_{1}=i_{2}\neq i_{3}=i_{4} \vee i_{1}=i_{3}\neq i_{2}=i_{4}
    \vee i_{1}=i_{4}\neq i_{2}= i_{3}
    \}$.

We discern an interesting structure in the original set $ I^{4}$ by partitioning it into four distinct subsets: $ I^{4}
=
\bigcup_{j=1}^{4}
 I_{j}^{4}$. Evidently, this decomposition is characterized by the fact that the sum of the cardinalities of these subsets, $\sum_{j=1}^{4}\left| I_{j}^{4}\right|$, yields $r+(4+3)r(r-1)+6r(r-1)(r-2)+r(r-1)(r-2)(r-3)=r^{3}$, which aligns with the cardinality of the original set $ I^{4}$. Moreover, it is apparent that $\bigcup_{j=1}^{4}
     I_{j}^{4}\subseteq I^{4}$, asserting the encompassing nature of the partition.

To firmly establish this conclusion, we resort to a proof by contradiction. Assume, for the sake of argument, that $  I^{4}
\neq
\bigcup_{j=1}^{4}
 I_{j}^{4}$. This implies the existence of an element $i\in  I^{4}$ that is not an element of $\bigcup_{j=1}^{4}
 I_{j}^{4}$. At this juncture, we arrive at a paradox: $\left| I^{4}\right|>\sum_{j=1}^{4}\left| I_{j}^{4}\right|$, contradicting the fundamental equality $\left| I^{4}\right|=\sum_{j=1}^{4}\left| I_{j}^{4}\right|$. Hence, we establish the validity of $ I^{4}
=
\bigcup_{j=1}^{4}
 I_{j}^{4}$ through contradiction.

Consequently, utilizing the definition of $\mathscr{P}$, this set equality, and the basic properties of summation, we can rewrite $4$-fold nested fuzzy summations as {\small
\begin{align*}
    &
    \sum_{i\in I^{4}}
    h_{i_{1}}h_{i_{2}}h_{i_{3}}h_{i_{4}}\Phi_{i_{1}i_{2}i_{3}i_{4}}
    \\
    &\quad
    =
    \frac{1}{4!}
     \sum_{i\in I^{4}}
     h_{i_{1}}h_{i_{2}}h_{i_{3}}h_{i_{4}}\mathscr{P}(\Phi_{i_{1}i_{2}i_{3}i_{4}})
      \\
    &\quad
    =
    \sum_{i\in I_{1}^{4}}
    h_{i_{1}}^{4}\Phi_{i_{1}i_{1}i_{1}i_{1}}
    +
    4\frac{3!}{4!}
    \sum_{i\in I_{2,1}^{4}}
    h_{i_{1}}^{3}h_{i_{2}}\mathscr{P}(\Phi_{i_{1}i_{1}i_{1}i_{2}})
    \\
    &\qquad
    +
    3\frac{2!2!}{4!}
    \sum_{i\in I_{2,2}^{4}}
    h_{i_{1}}^{2}h_{i_{2}}^{2}\mathscr{P}(\Phi_{i_{1}i_{1}i_{2}i_{2}})
    \\
    &\qquad
    +
    6\frac{2!}{4!}
    \sum_{i\in I_{3}^{4}}
    h_{i_{1}}^{2}h_{i_{2}}h_{i_{3}}\mathscr{P}(\Phi_{i_{1}i_{1}i_{2}i_{3}})
    \\
    &\qquad
    +
    \frac{1}{4!}
     \sum_{i\in I_{4}^{4}}
    h_{i_{1}}h_{i_{2}}h_{i_{3}}h_{i_{4}}\mathscr{P}(\Phi_{i_{1}i_{2}i_{3}i_{4}})
\end{align*}}
which results in \eqref{eq: quadruple sum}, where we used
{\small \begin{align*}
    &
    4
    \sum_{i\in I_{2,1}^{4}}
    h_{i_{1}}^{3}h_{i_{2}}\mathscr{P}(\Phi_{i_{1}i_{1}i_{1}i_{2}})
    \\
    &\quad
    =
    \sum_{i\in I_{2,1}^{4}}
    \left(
    h_{i_{1}}^{3}h_{i_{4}}\mathscr{P}(\Phi_{i_{1}i_{1}i_{1}i_{4}})
    +
    h_{i_{1}}^{3}h_{i_{3}}\mathscr{P}(\Phi_{i_{1}i_{1}i_{3}i_{1}})
    \right.
    \\
    &\qquad
    \left.
    +
    h_{i_{1}}^{3}h_{i_{2}}\mathscr{P}(\Phi_{i_{1}i_{2}i_{1}i_{1}})
    +
    h_{i_{1}}h_{i_{2}}^{3}\mathscr{P}(\Phi_{i_{1}i_{2}i_{2}i_{2}})
    \right)
    \\
    &
    3
     \sum_{i\in I_{2,2}^{4}}
    h_{i_{1}}^{2}h_{i_{2}}^{2}\mathscr{P}(\Phi_{i_{1}i_{1}i_{2}i_{2}})
    \\
    &\quad
    =
      \sum_{i\in I_{2,2}^{4}}
    \left(
    h_{i_{1}}^{2}h_{i_{3}}^{2}\mathscr{P}(\Phi_{i_{1}i_{1}i_{3}i_{3}})
    +
    h_{i_{1}}^{2}h_{i_{2}}^{2}\mathscr{P}(\Phi_{i_{1}i_{2}i_{1}i_{2}})
    \right.
    \\
    &\qquad
    \left.
    +
    h_{i_{1}}^{2}h_{i_{2}}^{2}\mathscr{P}(\Phi_{i_{1}i_{2}i_{2}i_{1}})
    \right)
    \\
    &6
    \sum_{i\in I_{4}^{3}}
    h_{i_{1}}^{2}h_{i_{2}}h_{i_{3}}\mathscr{P}(\Phi_{i_{1}i_{1}i_{2}i_{3}})
    \\
    &\quad
    =
     \sum_{i\in I_{4}^{3}}
    \left(
    h_{i_{1}}^{2}h_{i_{3}}h_{i_{4}}\mathscr{P}(\Phi_{i_{1}i_{1}i_{3}i_{4}})
    +
    h_{i_{1}}^{2}h_{i_{2}}h_{i_{4}}\mathscr{P}(\Phi_{i_{1}i_{2}i_{1}i_{4}})
    \right.
    \\
    &\qquad
    +
    h_{i_{1}}^{2}h_{i_{2}}h_{i_{3}}\mathscr{P}(\Phi_{i_{1}i_{2}i_{3}i_{1}})
    +
    h_{i_{2}}^{2}h_{i_{1}}h_{i_{4}}\mathscr{P}(\Phi_{i_{1}i_{2}i_{2}i_{4}})
    \\
    &\qquad
    \left.
    +
    h_{i_{2}}^{2}h_{i_{1}}h_{i_{3}}\mathscr{P}(\Phi_{i_{1}i_{2}i_{3}i_{2}})
    +
    h_{i_{3}}^{2}h_{i_{1}}h_{i_{2}}\mathscr{P}(\Phi_{i_{1}i_{2}i_{3}i_{3}})
    \right).
\end{align*}}
\section{Proof of Lemma \ref{lemma: multidimensional summation}}\label{app: proof of Lemma q-fold}
Utilizing set notation, we begin by representing the $q$-fold nested summations as
    $\sum_{i_{1}=1}^{r}
    \sum_{i_{2}=1}^{r}
    \cdots
    \sum_{i_{q}=1}^{r}
    h_{i}\Phi_{i}
    =
    \sum_{i\in I^{q}}h_{i}\Phi_{i}$, where $I^{q}$ represents the set of ordered $q$-tuples $i$ defined as $I^{q}=\{i\in\mathbb{N}_{r}^{q}\}$. We then proceed to partition of $I^{q}$ into $q$ distinct sets $I_{k}^{q}$, $k\in\mathbb{N}_{q}$ defined as the union over $\lambda\in\Lambda_{k}$ and over the collections $\{S_j\}_{j=1}^{q}$, where $S_j \subseteq \mathbb{N}_{q}$, $
    |S_{j}| = \lambda_{j}$, and $
    \bigcap_{j=1}^{q} S_j = \emptyset$.
 Specifically,
\begin{equation*}
\resizebox{1\hsize}{!}{$
     I_{k}^{q}
    =\bigcup_{\lambda\in\Lambda_{k}}
    \bigcup_{
    \substack{
    \{S_j\}_{j=1}^{q}, \\
    S_j \subseteq \mathbb{N}_{q}, \\
    |S_{j}| = \lambda_{j}, \\
    \bigcap_{j=1}^{q} S_j = \emptyset}
    }
    \left\{ i \;\middle|\;
    \forall j: \forall a,b \in S_j, \forall c \in \mathbb{N}_q \backslash \{a,b\}: i_a = i_b \neq i_c \right\}$}
\end{equation*}
for $k\in\mathbb{N}^{q}$.

This partitioning reveals an intriguing structure within the original set $I^{q}$: $I^{q}
=
\bigcup_{k=1}^{q}
I_{k}^{q}$. Drawing from Definition \ref{def:stirling_numbers} and Lemma \ref{lemma:stirling_powers}, this decomposition is characterized by the sum of the cardinalities of these subsets, $\sum_{k=1}^{q}\left|I_{k}^{q}\right|$, yields $s(q,1)\frac{r!}{(r-1)!}+s(q,2)\frac{r!}{(r-2)!}+\cdots+s(q,q
)\frac{r!}{(r-q)!}=r^{q}$, which coincides with the cardinality of the original set $I^{q}$. Furthermore, it is clear that $\bigcup_{k=1}^{q}
    I_{k}^{q}\subseteq I^{q}$, confirming the comprehensive nature of the partition.

To firmly establish this conclusion, we resort to a proof by contradiction. Assume, for the sake of argument, that $I^{q}
\neq
\bigcup_{k=1}^{q}
I_{k}^{q}$. This implies the existence of an element $i\in I^{q}$ that is not an element of $\bigcup_{k=1}^{q}
I_{k}^{q}$. At this juncture, we arrive at a paradox: $\left|I^{q}\right|>\sum_{k=1}^{q}\left|I_{k}^{q}\right|$, contradicting the fundamental equality $\left|I^{q}\right|=\sum_{k=1}^{q}\left|I_{k}^{q}\right|$. Hence, we establish the validity of $I^{q}
=
\bigcup_{k=1}^{q}
I_{k}^{q}$ through contradiction.

As a result, drawing upon this set equality, the definition of utilizing the definition of $\mathscr{P}$, and the general identities of summation, we can further express $q$-fold nested summation as
{\small
\begin{align*}
    \sum_{i\in I^{q}}
    h_{i}\Phi_{i}
    &=
     \sum_{\lambda\in\Lambda_{k},k\in\mathbb{N}_{q}}
    \binom{q}{\lambda_{1},\cdots,\lambda_{q}} \frac{1}{\mu(\lambda)!}
    \sum_{i_{1}=1}^{r}
    \sum_{\substack{
    i_{a},
    i_{a}\neq i_{b}\\
    a\in\mathbb{N}_{[2,k]}\\
    b\in\mathbb{N}_{a-1}}}
    h_{i}^{\lambda}\Phi_{i^{\lambda}}
   \\
    &
    =
     \sum_{\lambda\in\Lambda_{k},k\in\mathbb{N}_{q}}
    \binom{q}{\lambda_{1},\cdots,\lambda_{q}} \frac{1}{\mu(\lambda)!}
    \binom{q}{        \lambda_{1},\cdots,\lambda_{q}}^{-1}
    \\
    &\quad\times
      \sum_{i_{1}=1}^{r}
    \sum_{\substack{
    i_{a},
    i_{a}\neq i_{b}\\
    a\in\mathbb{N}_{[2,k]}\\
    b\in\mathbb{N}_{a-1}}}
    h_{i}^{\lambda}\mathscr{P}(\Phi_{i^{\lambda}}).
\end{align*}}
This completes the proof.
\section{Lemma \ref{lemma: multidimensional summation} and Theorem \ref{th: LMI q-dimensional} in the case of $4$-fold nested fuzzy summations}\label{app: example}
We begin by establishing the equivalence between the results of Lemma \ref{lemma: multidimensional summation} for the case when $q=4$ and Lemma \ref{lemma: quadruple summation}. According to Lemma \ref{lemma: multidimensional summation}, a four-dimensional summation can be represented as follows:
    $\sum_{i_{1}=1}^{r}
    \sum_{i_{2}=1}^{r}
    \sum_{i_{3}=1}^{r}
    \sum_{i_{4}=1}^{r}
    h_{i}\Phi_{i}
    = \sum_{\lambda\in\Lambda_{k},k\in\mathbb{N}_{4}}
    \frac{1}{\mu(\lambda)!}
      \sum_{i_{1}=1}^{r}
    \sum_{\substack{
    i_{a},
    i_{a}\neq i_{b}\\
    a\in\mathbb{N}_{[2,k]}\\
    b\in\mathbb{N}_{a-1}
    }
    }h_{i}^{\lambda}\mathscr{P}(\Phi_{i^{\lambda}})$.
To explore the right-hand side of the equation further, we introduce specific variable assignments derived from Definition \ref{def: the set of ordered $q$-tuple}. We denote the sets as $\Lambda_{1}=\{(4,0,0,0)\}$, $\Lambda_{2}=\{(3,1,0,0),(2,2,0,0)\}$, $\Lambda_{3}=\{(2,1,1,0)\}$, and $\Lambda_{4}=\{(1,1,1,1)\}$. Additionally, from Definition \ref{def:stirling_numbers}, we have the factorial values: $\mu((4,0,0,0))!=1$, $\mu((3,1,0,0))!=1$, $\mu((2,2,0,0))!=2$, $\mu((2,1,1,0))!=2$, and $\mu((1,1,1,1))!=24$. Utilizing these variables and Notations \ref{notation: multi-index} and \ref{notation: multi-summation}, we arrive at the expanded form:{\small
\begin{align*}
    &\sum_{\lambda\in\Lambda_{k},k\in\mathbb{N}_{4}}
    \frac{1}{\mu(\lambda)!}
      \sum_{i_{1}=1}^{r}
    \sum_{\substack{
    i_{a},
    i_{a}\neq i_{b}\\
    a\in\mathbb{N}_{[2,k]}\\
    b\in\mathbb{N}_{a-1}
    }
    }h_{i}^{\lambda}\mathscr{P}(\Phi_{i^{\lambda}})
    \\
    &= \sum_{i_{1}=1}^{r}
    h_{i_{1}}^{4}\mathscr{P}(\Phi_{i_{1}i_{1}i_{1}i_{1}})
    +
    \sum_{i_{1}=1}^{r}
    \sum_{\substack{
    i_{2},
    i_{2}\neq i_{1}
    }
    }h_{i_{1}}^{3}h_{i_{2}}\mathscr{P}(\Phi_{i_{1}i_{1}i_{1}i_{2}})
    \\&\quad+
    \frac{1}{2!}
   \sum_{i_{1}=1}^{r}
    \sum_{\substack{
    i_{2},
    i_{2}\neq i_{1}
    }
    }h_{i_{1}}^{2}h_{i_{2}}^{2}\mathscr{P}(\Phi_{i_{1}i_{1}i_{2}i_{2}})
    \\&\quad+
    \frac{1}{2!}
    \sum_{i_{1}=1}^{r}
    \sum_{\substack{
    i_{a},
    i_{a}\neq i_{b}\\
    a\in\mathbb{N}_{[2,3]}\\
    b\in\mathbb{N}_{a-1}
    }
    }h_{i_{1}}^{2}h_{i_{2}}h_{i_{3}}\mathscr{P}(\Phi_{i_{1}i_{2}i_{3}})
     \\&\quad+
    \frac{1}{4!}
    \sum_{i_{1}=1}^{r}
    \sum_{\substack{
    i_{a},
    i_{a}\neq i_{b}\\
    a\in\mathbb{N}_{[2,4]}\\
    b\in\mathbb{N}_{a-1}
    }
    }h_{i_{1}}h_{i_{2}}h_{i_{3}}h_{i_{4}}\mathscr{P}(\Phi_{i_{1}i_{2}i_{3}i_{4}})
\end{align*}}
which simplifies to the expression \eqref{eq: quadruple sum} presented in Lemma \ref{lemma: quadruple summation}. This concludes the demonstration of equivalence between the two lemmas.

Similarly, we can apply a similar approach to establish the equivalence between Theorem \ref{th: LMI q-dimensional} for the case of $q=4$ and Theorem \ref{th: LMI 4-dimensional}. By utilizing the sets \(\Lambda_{k}\), \(k\in\mathbb{N}_{4}\), and the factorials \(\mu(\lambda)!\), in conjunction with the conventions specified in Notations \ref{notation: multi-index} and \ref{notation: multi-summation}, we can observe that{\small
\begin{align*}
    &
    \eqref{eq: LMI q-dimensional}
    \text{ in Th. \ref{th: LMI q-dimensional}}
    \\
    &\quad\Leftrightarrow
        \Phi_{i_{1}i_{1}i_{1}i_{1}}
   +
    \sum_{
    \begin{smallmatrix}
    i_{2},
    i_{2}\neq i_{1}\\
    \end{smallmatrix}
    }
    \delta_{1,1,\sigma_{3}}
      \frac{1}{4}
    \left(
    3
    \left. \mathscr{P}(\Phi_{i_{1}i_{1}i_{1}i_{2}})\right|_{\begin{smallmatrix}
         i_{1}\mapsto i_{1}
         \\
         i_{1}\mapsto i_{1}
          \end{smallmatrix}}
    \right.
    \\
    &\qquad
    \left.
    +
    \left. \mathscr{P}(\Phi_{i_{1}i_{1}i_{1}i_{2}})\right|_{\begin{smallmatrix}
         i_{2}\mapsto i_{1}
         \\
         i_{1}\mapsto i_{2}
          \end{smallmatrix}}
    \right)
     +
   \frac{1}{2!}
    \sum_{
    \begin{smallmatrix}
    i_{2},
    i_{2}\neq i_{1}\\
    \end{smallmatrix}
    }
    \delta_{1,2,\sigma_{3}}
      \frac{1}{4}
      \\
      &\qquad\times
    \left(
    2
    \left. \mathscr{P}(\Phi_{i_{1}i_{1}i_{2}i_{2}})\right|_{\begin{smallmatrix}
         i_{1}\mapsto i_{1}
         \\
         i_{1}\mapsto i_{1}
          \end{smallmatrix}}
    +2
    \left. \mathscr{P}(\Phi_{i_{1}i_{1}i_{2}i_{2}})\right|_{\begin{smallmatrix}
         i_{2}\mapsto i_{1}
         \\
         i_{1}\mapsto i_{2}
          \end{smallmatrix}}
    \right)
    \\
    &\qquad
     +
   \frac{1}{2!}
    \sum_{
    \begin{smallmatrix}
    i_{a},
    i_{a}\neq i_{b}\\
    a\in\mathbb{N}_{[2,3]}\\
    b\in\mathbb{N}_{a-1}
    \end{smallmatrix}
    }
    \delta_{2,1,\sigma_{3}}
      \frac{1}{4}
      \left(
    2\left. \mathscr{P}(\Phi_{i_{1}i_{1}i_{2}i_{3}})\right|_{\begin{smallmatrix}
         i_{1}\mapsto i_{1}
         \\
         i_{1}\mapsto i_{1}
          \end{smallmatrix}}
    \right.
    \\
    &\qquad
    +
    \left.
    \left. \mathscr{P}(\Phi_{i_{1}i_{1}i_{2}i_{3}})\right|_{\begin{smallmatrix}
         i_{2}\mapsto i_{1}
         \\
         i_{1}\mapsto i_{2}
          \end{smallmatrix}}
    +\left. \mathscr{P}(\Phi_{i_{1}i_{1}i_{2}i_{3}})\right|_{\begin{smallmatrix}
         i_{3}\mapsto i_{1}
         \\
         i_{1}\mapsto i_{3}
          \end{smallmatrix}}
    \right)
    \\
   &\qquad  +
   \frac{1}{4!}
    \sum_{
    \begin{smallmatrix}
    i_{a},
    i_{a}\neq i_{b}\\
    a\in\mathbb{N}_{[2,4]}\\
    b\in\mathbb{N}_{a-1}
    \end{smallmatrix}
    }
    \delta_{3,1,\sigma_{3}}
      \frac{1}{4}
      \left(
    \left. \mathscr{P}(\Phi_{i_{1}i_{2}i_{3}i_{4}})\right|_{\begin{smallmatrix}
         i_{1}\mapsto i_{1}
         \\
         i_{1}\mapsto i_{1}
          \end{smallmatrix}}
          \right.
    \\
    &\qquad
    \left.
    +
    \left. \mathscr{P}(\Phi_{i_{1}i_{2}i_{3}i_{4}})\right|_{\begin{smallmatrix}
         i_{2}\mapsto i_{1}
         \\
         i_{1}\mapsto i_{2}
          \end{smallmatrix}}
          \right.
          \left.
    +\left. \mathscr{P}(\Phi_{i_{1}i_{2}i_{3}i_{4}})\right|_{\begin{smallmatrix}
         i_{3}\mapsto i_{1}
         \\
         i_{1}\mapsto i_{3}
          \end{smallmatrix}}
    \right.
    \\
    &\qquad
    \left.
    +\left. \mathscr{P}(\Phi_{i_{1}i_{2}i_{3}i_{4}})\right|_{\begin{smallmatrix}
         i_{4}\mapsto i_{1}
         \\
         i_{1}\mapsto i_{4}
          \end{smallmatrix}}
    \right)
    \prec0
    \\
    &\quad
    \Leftrightarrow
    \eqref{eq: LMI quadruple}
    \text{ in Th. \ref{th: LMI 4-dimensional}}
\end{align*}}
where we used $\sigma_{1}(2)=1$,
    $\sigma_{1}(3)=2$,
    $\sigma_{1}(4)=3$,
    $\sigma_{2}(2,(3,1,0,0)))=1$,
    $\sigma_{2}(2,(2,1,0,0)))=2$,
    $\sigma_{2}(3,\lambda)=3$,
    $\sigma_{2}(4,\lambda)=4$,
    $\sigma_{3}(2,i)\in\{1,2,\cdots,r-1\}$,
    $\sigma_{3}(3,i)\in\{1,2,\cdots,(r-1)(r-2)\}$, and
    $\sigma_{3}(4,i)\in\{1,2,\cdots,(r-1)(r-2)(r-3)\}$
in Definition \ref{def: general sigma function}, $\mathscr{P}(\Phi_{i_{1}i_{1}i_{2}i_{2}})=\mathscr{P}(\Phi_{i_{2}i_{2}i_{1}i_{1}})$, and
$ \mathscr{P}(\Phi_{i_{1}i_{2}i_{3}i_{4}})
    =\mathscr{P}(\Phi_{i_{2}i_{1}i_{3}i_{4}})
    =\mathscr{P}(\Phi_{i_{3}i_{2}i_{1}i_{4}})
    =\mathscr{P}(\Phi_{i_{4}i_{2}i_{3}i_{1}})$. This concludes the demonstration of equivalence between the two theorems.
\bibliographystyle{IEEEtran}
\bibliography{reference}

\end{document}